\documentclass[MSNbibl,number,citesort,seceqn,dvips]{arxbj}
\usepackage{upgreek}
\usepackage{graphicx}

\aid{0}
\volume{18}
\issue{3}
\pubyear{2012}
\firstpage{747}
\lastpage{763}
\doi{10.3150/12-BEJ346} %kopijuoti is PTS

\makeatletter
\newproclaim{asum}{Assumption}
\newtheorem{Prop}{Proposition}[section]
\newtheorem{Corol}{Corollary}[section]
\newcommand{\R}{\mathbb R}

\newcommand{\pr}{^{\prime}}

\newcommand{\Xb}{\mathbf{X}}
\newcommand{\Yb}{\mathbf{Y}}
\newcommand{\Zb}{\mathbf{Z}}
\newcommand{\Ob}{{\mathbf{O}}}
\newcommand{\Mb}{\mathbf{M}}
\newcommand{\varthetab}{{\boldsymbol\vartheta}}
\newcommand{\omegab}{{\boldsymbol\omega}}
\newcommand{\Sigmab}{\boldsymbol{{\Sigma}}}
\newcommand{\deltab}{{\boldsymbol\delta}}
\newcommand{\psib}{{\boldsymbol\psi}}
\newcommand{\mub}{{\boldsymbol\mu}}
\newcommand{\vecb}{\operatorname{vec}}
\newcommand{\vechb}{\operatorname{vech}}
\newcommand{\zerob}{{\mathbf0}}
\newcommand{\varphib}{{\boldsymbol\varphi}}
\makeatother

\begin{document}
\begin{frontmatter}

\title{Skew-symmetric distributions and Fisher information --
a tale of two densities}
\runtitle{Skew-symmetric distributions and Fisher information -- a tale
of two densities}

\begin{aug}
%%%% inicialai - be tarpu
\author[1]{\fnms{Marc} \snm{Hallin}\corref{}\thanksref{1}\ead[label=e1]{mhallin@ulb.ac.be}} \and
\author[2]{\fnms{Christophe} \snm{Ley}\thanksref{1,2}\ead[label=e2]{chrisley@ulb.ac.be}}
%%\runauthor{} %% auto
\address[1]{E.C.A.R.E.S., CP 114, Universit\'{e} Libre de Bruxelles, 50
Avenue F.D. Roosevelt, 1050 Brussels, Belgium. \printead{e1}}
\address[2]{D\' epartement de Math\' ematique, CP 210, Universit\'{e}
Libre de Bruxelles, Boulevard du Triomphe, 1050 Brussels, Belgium.
\printead{e2}}
\end{aug}

% HISTORY:
\received{\smonth{7} \syear{2010}}
\revised{\smonth{11} \syear{2010}}

% ABSTRACT
%
\begin{abstract}
Skew-symmetric densities recently received much attention in the literature,
giving rise to increasingly general families of univariate and multivariate
skewed densities. Most of those families, however, suffer from the inferential
drawback of a potentially singular Fisher information in the vicinity
of symmetry.
All existing results indicate that Gaussian densities (possibly after
restriction to
some linear subspace) play a special and somewhat intriguing role in
that context.
We dispel that widespread opinion by providing a full characterization,
in a general
multivariate context, of the information singularity phenomenon,
highlighting its relation
to a possible link between symmetric kernels and skewing functions -- a
link that can be
interpreted as the mismatch of two densities.
\end{abstract}

% KEYWORDS
%

\end{frontmatter}

%s1 ###
\section{Introduction}
Models for skewed distributions have become increasingly popular in
recent years, as they provide a much better fit for data presenting
some departure from normality, and from symmetry in general. Many of
the proposed models in the literature allow for a~continuous variation
from symmetry to asymmetry, regulated by some finite-dimensional parameter.

The success of those skewed distributions started with the seminal
papers by Azzalini \cite{A85,A86} introducing the scalar \emph
{skew-normal} model, which embeds the univariate normal distributions
into a flexible parametric class of (possibly) skewed distributions.
More formally, a random variable $X$ is said to be skew-normal with
location parameter $\mu\in\R$, scale parameter $\sigma\in\R^+_0$ and
skewness parameter $\delta\in\R$ if it admits the probability density
function (p.d.f.)
%
%e1.1 ###
\begin{equation}\label{HLSN}
x\mapsto
%f_{\mu,\delta}(x)=
2\sigma^{-1}\phi\bigl(\sigma^{-1}(x-\mu)\bigr) \Phi\bigl(\delta\sigma
^{-1}(x-\mu)\bigr),\qquad
x\in\R,
\end{equation}
where $\phi$ and $\Phi$ respectively denote the p.d.f. and cumulative distribution function (c.d.f.) of a standard
normal distribution. Besides their many appealing features, however,
skew-normal densities unfortunately also suffer from an unpleasant
inferential drawback: in the vicinity of symmetry, that is, at $\delta
=0$, the Fisher information matrix for the three-parameter density
(\ref
{HLSN}) is singular -- typically, with rank 2 instead of 3.
Consequently, skew-normal distributions happen to be problematic from
an inferential point of view, since that singularity violates the
assumptions for standard Gaussian asymptotics and precludes, at first
sight, any nontrivial test of the null hypothesis of symmetry. Such a~situation has been studied by Rotnitzky \textit{et al.} \cite{Ro00}, who show
that one of the parameters then cannot be estimated at the usual
root-$n$ rate, while the limit distribution of maximum likelihood
estimators might be bimodal.

This Fisher singularity problem, however, did not hamper the success of
skew-normal densities among practitioners, while theoretical extensions
were developing into various directions. Azzalini and Dalla Valle
\cite{AD96} and Azzalini and Capitanio \cite{AC99} consider multivariate
skew-normal distributions resulting from replacing in (\ref{HLSN}) the
univariate normal kernel $\phi$ with its $k$-variate version $\phi_k$.
In the same paper, Azzalini and Capitanio also propose substituting an
elliptical kernel $f_k$ for the normal one $\phi_k$, and replacing the
skewing factor $\Phi$ in (\ref{HLSN}) with an arbitrary, possibly
non-Gaussian, univariate symmetric c.d.f.~$G_1$. The resulting
distributions are called \textit{skew-elliptical}. The class of
skew-elliptical distributions is also studied in detail by Branco and
Dey \cite{BD01}, based, however, on a slightly different definition. Genton
and Loperfido \cite{GL05} introduce a concept of \textit{generalized
skew-elliptical} distributions encompassing all previous ones, where
arbitrary skewing functions (not necessarily c.d.f.s, but satisfying
c.d.f.-type conditions) can be used in conjunction with the elliptical
kernel $f_k$. Finally, Azzalini and Capitanio \cite{AC03} (who also propose
the nowadays commonly adopted definition of multivariate skew-$t$
distributions), then Wang \textit{et al.} \cite{WBG04} are relaxing the
assumption of elliptically symmetric kernels into a weaker assumption
of central symmetry, leading to multivariate \emph{skew-symmetric}
densities of the form
%
%e1.2 ###
\begin{eqnarray}\label{HLSS}
{\mathbf{x}}&\mapsto& f_\varthetab^\Pi(
{\mathbf{x}})=f_{\mub,\Sigmab,\deltab}^\Pi({\mathbf{x}})\nonumber
\\[-8pt]
\\[-8pt]
&:=&2
|\Sigmab|^{-1/2}f\bigl(\Sigmab^{-1/2}({\mathbf{x}}-\mub)\bigr)\Pi
\bigl(\Sigmab^{-1/2}({\mathbf{x}}-\mub
),\deltab\bigr), \qquad{\mathbf{x}}\in\R^k,\nonumber
\end{eqnarray}
where
\begin{enumerate}[(b)]
\item[(a)] $\mub\in\R^k$ is a location parameter, $\Sigmab\in
\mathcal
{S}_k$ (throughout, $|\Mb|$ denotes the determinant and $\Mb^{1/2}$ the
symmetric square-root of any $\Mb$ in the class $\mathcal{S}_k$ of
symmetric positive definite $k\times k$ matrices) a scatter matrix,
while $\deltab\in\R^k$ plays the role of a skewness parameter;
\item[(b)] the \emph{symmetric kernel} $f$ is a centrally symmetric
nonvanishing p.d.f., meaning that $0\neq f(-{\mathbf
{z}})=f({\mathbf{z}})$, $ {\mathbf{z}}\in\R
^k$, and
\item[(c)] the \emph{skewing function} $\Pi\dvtx\R^k\times\R
^k\rightarrow
[0,1]$ satisfies $\Pi(-{\mathbf{z}},\deltab)+\Pi
({\mathbf{z}},\deltab)=1$, ${\mathbf
{z}},\deltab
\in\R^k$, and $\Pi({\mathbf{z}},{\boldsymbol0})=1/2$,
$ {\mathbf{z}}\in\R^k$.
\end{enumerate}

This definition is the one we are adopting in the sequel. While $\Pi
({\mathbf{z}}
,\deltab)$, in most practical situations, is of the simple form $\Pi
(\deltab^{\prime}{\mathbf{z}})$, with $\Pi\dvtx \R
\rightarrow[0,1]$,
Wang \textit{et al.} \cite{WBG04} actually do not consider any specific
$\deltab$-parameterization. Our parametric approach (with the
regularity assumptions \hyperref[as2]{(A2)}--\hyperref[as2plus]{(A2$^+$)} and \hyperref[asb2]{(B2)}--\hyperref[asb2plus]{(B2$^+$)} of Sections~\ref{univtale}
and~\ref{multitale}, resp.) is in the spirit of -- if not at the same level
of mathematical generality as -- the differentiable path and tangent
space approach taken in the local and asymptotic treatment of
semiparametric models (see, e.g., Chapter~25 of van der Vaart \cite{vdV00}).
Also, the condition that $f$ is a nonvanishing density is not imposed
by Wang \textit{et al.} \cite{WBG04}; we are adding that requirement in order
to avoid inessential complications related with bounded and
parameter-dependent supports. For further information about
skew-symmetric models and related topics, we refer the reader to the
recent monograph by Genton \cite{G04}, and to the review papers Arnold and
Beaver \cite{AB02} and Azzalini \cite{A05}.

The issue of singular Fisher information runs like a red thread through
all those developments.
Mentioned, from the very beginning, in Azzalini \cite{A85} itself, it is
discussed, in the univariate and multivariate skew-normal context, by
Azzalini and Capitanio \cite{AC99}, Pewsey~\cite{P00}, Chiogna \cite{C05} and
Arellano-Valle and Azzalini~\cite{AA08}.
The same issue has been considered in various subclasses of
skew-symmetric distributions. Pewsey \cite{P06} and Azzalini and Genton
\cite{AG08} establish that the singularity problem remains after replacement
of the c.d.f.~$\Phi$ in (\ref{HLSN}) with any c.d.f. $H$ satisfying mild
regularity assumptions. DiCiccio and Monti~\cite{DM04} prove that, within
the class of univariate skew-exponential power distributions of
Azzalini \cite{A86}, the normal kernels are the only ones suffering from
singular Fisher information. The same result is shown to hold true for
two classes of scalar skew-$t$ distributions by G\'{o}mez \textit{et al.}
\cite{GTB07} and DiCiccio and Monti \cite{DM10}.
The multivariate counterparts of these statements are provided in Ley
and Paindaveine \cite{LP10a,LP10b}, respectively.

Finally, the very general (still a special case of (\ref{HLSS}),
though) class of multivariate skew-symmetric densities of the form
%
%e1.3 ###
\begin{equation}\label{leypaindav}
{\mathbf{x}}\mapsto2|\Sigmab|^{-1/2}f\bigl(\Sigmab
^{-1/2}({\mathbf{x}}-\mub)\bigr)\Pi\bigl(\deltab
^{\prime}\Sigmab^{-1/2}({\mathbf{x}}-\mub)\bigr),\qquad
{\mathbf{x}}\in\R^k,
\end{equation}
encompassing all previous cases, is considered in Ley and Paindaveine
\cite{LP10a}, who characterize, for each possible value $1\leq m\leq k$ of
the Fisher information rank deficiency, the form of the symmetric
kernels giving rise to such deficiency. Here again, Gaussian kernels
are playing a very special role. In the univariate setup and within the
subclass of multivariate generalized skew-elliptical distributions,
only the skew-normal densities are affected by the singularity problem.
Although results in the fully general (for densities of the form (\ref
{leypaindav})) multivariate case are more complex, only kernels
exhibiting Gaussian restrictions on some $m$-dimensional linear
subspaces can lead to degenerate Fisher information.

A tentative remedy to that singularity problem was suggested by
Azzalini himself who, as early as 1985, in his original paper, proposes
a reparametrization of skew-normal families, the so-called \textit
{centered parametrization}, under which Fisher information matrices
remain full-rank. The multivariate version of that reparametrization is
examined in detail by Arellano-Valle and Azzalini \cite{AA08}. That
solution, however, never really caught up in practice, partly because
the structure of the skewing mechanism, hence of the resulting
skew-normal family, under the new parametrization, loses much of its
simplicity (certainly so in the multivariate context), partly because
of its limitation to skew-normal families. Azzalini and Genton \cite{AG08}
therefore once again emphasize the need for a clarification of the
Fisher singularity phenomenon in order to ``\textit{remove, or at least
alleviate, the necessity of an alternative parametrization.}''

The objective of the present paper is to provide such a clarification.
While all comments and existing results, in this singular Fisher
information issue, seemed to be pointing at some special status for
normal kernels and, consequently, skew-normal distributions, we
completely dispel the idea of any particular role of Gaussian kernels.
Turning to the fully general class of skew-symmetric densities
described in (\ref{HLSS}), we show indeed that information deficiency
actually originates in an unfortunate mismatch between $f$ and $\Pi
$ -- more specifically, between two densities, the kernel $f$ and an
exponential density $g_\Pi$ associated with the skewing function $\Pi$.

A tale of two densities, thus, rather than a Gaussian mystery\ldots\

The paper is organized as follows. Section~\ref{univtale} deals with
the univariate setup, where the singularity problem is simple, as the
rank of the three-parameter Fisher information matrix only can be 3 or
2. The result is derived in an informal way, and some examples of
skewing functions are treated in Section~\ref{univexos}. A more formal
statement of the general solution is provided for the multivariate
setup in Section~\ref{multitale}, along with some examples in
Section~\ref{multiexos}. Final comments and conclusions are given in
Section~\ref{fcomments}.
%
%s2 ###
\section{The univariate setup}\label{univar}
%
%s2.1 ###
\subsection{A tale of two densities \ldots}\label{univtale} We start
by analyzing the information singularity problem in the univariate
case. To do so, consider the class of skew-symmetric
probability density families of the form
%
%e2.1 ###
\begin{equation}
x\mapsto f_{\varthetab}^\Pi(x)= f_{\mu,\sigma,\delta}^\Pi(x) :=2
\sigma
^{-1}f\bigl(\sigma^{-1}(x-\mu)\bigr)\Pi\bigl(\sigma^{-1}(x-\mu),\delta\bigr),\qquad x\in\R
,\label{skew1}
\end{equation}
with $\varthetab:= (\mu,\sigma,\delta)\pr$, where $\mu\in\R$ is a
location parameter, $\sigma\in\R^+_0$ a scale parameter and $\delta
\in\R
$ an asymmetry parameter.

The symmetric kernel $f\dvtx \R\rightarrow\R^+$ in (\ref{skew1}) is a
nonvanishing symmetric \textit{standardized} p.d.f., that is, a
probability density function such that $f(z)=f(-z)\neq0 $ for all
$z\in
\R$, with scale parameter one -- an identification constraint for
$\sigma
$ that does not imply any loss of generality. Classical
standardization, with a constraint of the form $\int_{-\infty}^\infty
z^2f(z)\,\mathrm{d}z = 1$, involves the variance of $Z$ with p.d.f. $f$; the scale
parameter $\sigma^2$ then is the mean squared deviation $\mathrm{{E}}
[(X-\mu)^2]$ with respect to $\mu$ of $X$ with p.d.f. $f_{\mu,\sigma
,0}^\Pi$. If moment assumptions are to be avoided, one may rather
consider, for instance, medians of squares, with an identification
constraint of the form $\int_{-\infty}^1f(z)\,\mathrm{d}z=0.75$: if $X$ has p.d.f.
$f_{\mu,\sigma,0}^\Pi$, $\sigma$ then is the median of the absolute
deviation $\vert X-\mu\vert$, which exists irrespective of the density
of $X$. Other quantiles of $\vert X-\mu\vert$ would enjoy similar properties.
%Except for possible existence constraints it may place on $f$, the
%choice of this identification constraint as no impact on the results
%of this paper. We therefore
We throughout assume that such an identification constraint, hence a
concept of scale, has been adopted. That choice, however, is completely
arbitrary, and any element in the scale family of p.d.f.'s of the form
(\ref{skew1}) with $\mu=\delta= 0$ could be chosen as the reference density
characterizing unit scale -- hence could serve as a symmetric kernel for
the same skew-symmetric family. As we shall see, that choice has no
impact on the results of this paper.

The second factor in (\ref{skew1}) is a skewing function, namely, a
function $\Pi\dvtx\R\times\R\rightarrow[0,1]$
% \textit{a.e.} continuous, differentiable in its second argument at 0
%and
such that $\Pi(-z,\delta)+\Pi(z,\delta)=1$ for all $z,\delta\in\R$,
and $\Pi(z,0)=1/2$ for all $z\in\R$. Traditional choices involve
$\Pi
(z,\delta)=\Phi(\delta z)$ (skew-normal distributions, Azzalini \cite{A85}),
$\Pi(z,\delta)=\Phi(\delta\operatorname{sign}(z)|z|^{\alpha
/2}(2/\alpha)^{1/2})$ (skew-exponential power distributions, Azzalini
\cite{A86}) or $\Pi(z,\delta)=G(\delta z)$ for any symmetric univariate c.d.f.
$G$ (skew-symmetric distributions, Azzalini and Capitanio \cite{AC99}). The
class of skewing functions considered here is much broader.
%: except for the regularity assumptions below, we do not impose any
%restriction on $\Pi$.

The regularity assumptions we are making on $f$ and $\Pi$ are as follows.

\begin{asum*}[(A1)]\label{as1} The mapping $z\mapsto f(z)$ is differentiable,
with derivative $\dot f$ such that, letting $\varphi_f:=-\dot{f}/f$,
the information quantity for location $\sigma^{-2}
\mathcal{I}_f$, with
\[
\mathcal{I}_f:=\int_{-\infty}^{\infty}\varphi^2_f(z) f(z)\,\mathrm{d} z,
\]
is finite.
\end{asum*}
\begin{asum*}[(A1$^+$)]\label{as1plus} Same as \hyperref[as1]{(\textup{A}1)}, but the information quantity
for scale $\sigma^{-2}\mathcal{J}_f$, with
\[
\mathcal{J}_f:= \int_{-\infty}^{\infty}\bigl(z\varphi_f(z)-1\bigr)^2 f(z)\,\mathrm{d} z,
\]
moreover is finite.
\end{asum*}
\begin{asum*}[(A2)]\label{as2} \emph{(i)}  The mapping $z\mapsto\Pi(z,\delta)$ is
differentiable, and its derivative equals 0 at $\delta=0$. \emph{(ii)} The
mapping $\delta\mapsto\Pi(z,\delta)$ is differentiable at $\delta= 0$
for all $z\in\R$, with derivative (at $\delta= 0$) ${{\partial}_{
\delta}\Pi(z,\delta)|_{\delta=0}=:\psi(z)}$ such that $z\mapsto
\psi(z)$
admits a primitive, denoted as $\Psi$.
\end{asum*}
\begin{asum*}[(A2$^+$)]\label{as2plus} Same as \hyperref[as2]{(\textup{A}2)}, but the quantity
\[
\int_{-\infty}^\infty\psi^2(z)f(z)\,\mathrm{d}z
\]
moreover is finite.
\end{asum*}

These assumptions essentially guarantee the existence and finiteness of
Fisher information at $\delta=0$; the differentiability and
integrability conditions could be relaxed into weaker differentiability
properties such as quadratic mean differentiability. This small gain of
generality, however, would require a generalized definition of
information (in the Le Cam style), with non-negligible technical
complications. For the sake of simplicity, we stick to a more
traditional approach and the traditional definition of Fisher
information. Note that this definition
%of a \textit{population information matrix} (or \textit{expected}
%information matrix)
differs from the one, used by some authors, of an \textit{observed
Fisher information}, that is, the empirical value of the matrix of
negative second-order derivatives of the log-likelihood evaluated at
the maximum likelihood estimator of the parameters.

Under Assumptions~\hyperref[as1]{(A1)} and~\hyperref[as2]{(A2)},
the \emph{score vector} ${\boldsymbol\ell}_{f;\varthetab}$, at $
(\mu, \sigma,
0)^\prime=:\varthetab_0$, takes the form
\begin{eqnarray*}
{\boldsymbol\ell}_{f;\varthetab_0}(x)&:=&
\operatorname{grad}_\varthetab\log f_{\varthetab}^\Pi(x) \vert
_{\varthetab_0}=: (
\ell_{f;\varthetab_0}^{1}(x),
\ell_{f;\varthetab_0}^{2}(x),
\ell_{f;\varthetab_0}^{3}(x)
)\pr\\
&=& \pmatrix{
\sigma^{-1}\varphi_f\bigl(\sigma^{-1}(x-\mu)\bigr)\cr
\sigma^{-1}\bigl(\sigma^{-1}(x-\mu)\varphi_f\bigl(\sigma^{-1}(x-\mu)\bigr)-1\bigr)\cr
2\psi\bigl(\sigma^{-1}(x-\mu)\bigr)}
,
\end{eqnarray*}
where the factor 2 in $\ell_{f;\varthetab_0}^{3}$ follows from the
fact that $\Pi(z,0)=1/2$ for all $z\in\R$. Assumption~\hyperref[as2]{(A2)}(i) is a mild
requirement which, in regular models, readily follows from the fact
that $\Pi(z,0)=1/2$, and ensures that the skewing function $\Pi$ plays
no role in the score functions for $\mu$ and $\sigma$ at $\delta=0$.

Under Assumptions~\hyperref[as1plus]{(A1$^+$)} and \hyperref[as2plus]{(A2$^+$)},
the $3\times3$ Fisher information matrix for $(\mu, \sigma, \delta)$
exists, and takes the form
\[
\boldsymbol{{\Gamma}}_{f;\varthetab_0}:=\sigma^{-1} \int_{-\infty
}^{\infty} {\boldsymbol\ell}_{f;\varthetab_0}(x){\boldsymbol\ell
}^{\prime
}_{f;\varthetab_0}(x) f\bigl(\sigma^{-1} (x-\mu)\bigr)\,\mathrm{d}x
=:
\pmatrix{
\gamma_{f;\varthetab_0}^{11}&0&\gamma_{f;\varthetab_0}^{13}\cr
0&\gamma_{f;\varthetab_0}^{22}&0\cr
\gamma_{f;\varthetab_0}^{13}&0&\gamma_{f;\varthetab_0}^{33}
}
,
\]
with
\[
\gamma_{f;\varthetab_0}^{11}=\sigma^{-2}
\mathcal{I}_f ,\qquad
%%\int_\R\varphi_f^2(x)f_{0,1,0}^\Pi(x),
\gamma_{f;\varthetab_0}^{22}=\sigma^{-2}\mathcal{J}_f,\qquad
\gamma_{f;\varthetab_0}^{33}=4\int_{-\infty}^\infty\psi^2(z)f(z)\,\mathrm{d} z
\]
and
\[
\gamma_{f;\varthetab_0}^{13}=2\sigma^{-1} \int_{-\infty}^\infty
\varphi
_f(z)\psi(z)f(z)\,\mathrm{d} z.
\]
The zeroes in $\boldsymbol{{\Gamma}}_{f;\varthetab_0}$ are easily
obtained by noting that $\ell^1_{f;\varthetab_0}$ and $\ell
^3_{f;\varthetab_0}$ are antisymmetric functions of $(x-\mu)$, whereas
$\ell^2_{f;\varthetab_0}$ is symmetric with respect to the same quantity.

It then trivially follows that singularity of $\boldsymbol{{\Gamma
}}_{f;\varthetab_0}$
only can be due to the singularity of the $2\times2$ submatrix
%cannot be due to the scale-part, hence we can, without any loss of
%generality, restrict our attention to the submatrix
%
\[
\boldsymbol{{\Gamma}}_{f;\varthetab_0}^{0}:= \pmatrix{
\gamma_{f;\varthetab_0}^{11}&\gamma_{f;\varthetab_0}^{13}\cr
\gamma_{f;\varthetab_0}^{13}&\gamma_{f;\varthetab_0}^{33}
},
\]
the existence of which, however, only requires Assumptions~\hyperref[as1]{(A1)} and
\hyperref[as2plus]{(A2$^+$)}. Clearly,
%is singular iff its determinant equals zero (note that ${\boldsymbol
either $\boldsymbol{{\Gamma}}_{f;\varthetab_0}^{0}$ is full-rank or,
in case $\gamma_{f;\varthetab_0}^{11}\gamma_{f;\varthetab
_0}^{33}=(\gamma_{f;\varthetab_0}^{13})^2$, it has rank 1.
%The latter condition is equivalent to .

Now, the Cauchy--Schwarz inequality implies that $(\gamma_{f;\varthetab
_0}^{13})^2\leq\gamma_{f;\varthetab_0}^{11}\gamma_{f;\varthetab
_0}^{33}$, with equality if and only if
%
%e2.2 ###
\begin{equation}\label{phipsi}\varphi_f=a \psi \qquad f\mbox{-a.s.
(equivalently, Lebesgue-a.e.)}
\end{equation}
for some constant $a\in\R$. It thus follows that ${\boldsymbol\Gamma
}_{f;\varthetab_0}^0$ is singular for any $ \varthetab_0 = (\mu,
\sigma, 0)^\prime$ if and only if (\ref{phipsi}) is satisfied
%$\varphi_f=a \psi$ Lebesgue-a.e.
for some $a\in\R$. This holds under Assumptions \hyperref[as1]{(A1)} and~\hyperref[as2plus]{(A2$^+$)}. If
Assumption~\hyperref[as1]{(A1)} is reinforced into \hyperref[as1plus]{(A1$^+$)}, the $2\times2$
singularity of ${\boldsymbol\Gamma}_{f;\varthetab_0}^0$ in turn is equivalent
to the $3\times3$ singularity of $\boldsymbol{{\Gamma
}}_{f;\varthetab_0}$.
Replacing $\varphi_f$ with its definition, the necessary and
sufficient condition
$\varphi_f=a \psi$
yields
a first-order differential equation whose solutions are of the form
$f(x)=c \exp(-a\Psi(x))$ for some $a\in\R$, where $\Psi$ is a
primitive of~$\psi$ and $c\in\R^+$ an integration constant.

Summing up, let the couple $(f,\Pi)$ satisfy Assumptions~\hyperref[as1plus]{(A1$^+$)} and~\hyperref[as2plus]{(A2$^+$)}:
$\boldsymbol{{\Gamma}}_{f;\varthetab_0}$ is singular for all
$\varthetab_0$ if and only if the symmetric kernel $f$ belongs to the
\textit{exponential family}
%
%e2.3 ###
\begin{equation}\label{expfam}
\mathcal{E}_\Psi:=
\biggl\{
g_a:=\exp(-a\Psi)\Big/\int_{-\infty}^\infty\exp(-a\Psi(z))\,\mathrm{d}z \Bigm| a \in\mathcal{A}
\biggr\}
\end{equation}
with \textit{minimal sufficient statistic} $\Psi$, \textit{natural
parameter} $-a$, and \textit{natural parameter space}
\[
\mathcal{A}:= \biggl\{ a\in\R\mbox{ such that } \int_{-\infty}^\infty
\exp
(-a\Psi(z))\,\mathrm{d}z<\infty
\biggr\}.
\]
The same statement can be made under Assumptions~\hyperref[as1]{(A1)}
and~\hyperref[as2plus]{(A2$^+$)}
about the singularity of~$\boldsymbol{{\Gamma}}_{f;\varthetab_0}^0$.

Note that $\mathcal{A}$, as the natural parameter space of an exponential
family, is an open interval of $\R$. The unique value $a_\Pi$ of
$a\in
\mathcal{A}$ such that $f$ and $g_{a_\Pi}$ coincide, if any, is entirely
determined by the standardization constraint on $f$. If the classical
variance-based standardization is adopted, then $a_\Pi$ is solution of
the equation
\[
\int_{-\infty}^\infty z^2 \exp(-a\Psi(z))\,\mathrm{d}z=\int_{-\infty}^\infty
\exp(-a\Psi(z))\,\mathrm{d}z.
\]
If standardization is imposed via medians of squares, $a_\Pi$ is
solution of
\[
\int_{-\infty}^1 \exp(-a\Psi(z))\,\mathrm{d}z= 3\int^{\infty}_1 \exp(-a\Psi
(z))\,\mathrm{d}z .
\]

Letting $f_\sigma(x):=\sigma^{-1}f(x/\sigma)$, $\sigma\in\R^+_0$, also
note that $f\in\mathcal{E}_\Psi$ if and only if $f_\sigma\in
\mathcal{E}_{\Psi
{ \circ}\sigma^{-1}}$, where $\mathcal{E}_{\Psi{\circ
}\sigma
^{-1}}$ stands for the exponential family with minimal sufficient
statistic $\Psi{\circ}\sigma^{-1}\dvtx z\mapsto\Psi
(\sigma^{-1}z)$. It is easy to see that both conditions moreover
determine the same~$a_\Pi$, which confirms that the arbitrary choice of
a scale parameter has no impact on the result.

As a consequence of those results, it follows that, for any symmetric
density $f$ satisfying Assumption~\hyperref[as1plus]{(A1$^+$)} (resp., Assumption~\hyperref[as1]{(A1)}),
there exists a skewing function $\Pi_f$ (infinitely many of them,
actually) such that $\boldsymbol{{\Gamma}}_{f;\varthetab_0}$ (resp.,
$\boldsymbol{{\Gamma}}_{f;\varthetab_0}^0$) exists and is singular for
any $\varthetab_0$; among them, with $a_\uppi= \sqrt{2\uppi}$, $\Pi_f
(z,\delta):=\Phi(\delta\varphi_f(z))$, for which
Assumption~\hyperref[as2plus]{(A2$^+$)} holds.

The converse is slightly more subtle. Let $\Pi$ be a skewing function
satisfying Assumption~\hyperref[as2]{(A2)}; a function $\Psi$ with derivative $\psi$
thus exists, which automatically satisfies $\Psi(z)=\Psi(-z)$. If there
exists a density $g_a$ in the corresponding exponential family (\ref
{expfam}) such that $\int_{-\infty}^{\infty}\psi^2(z)g_a(z)\,\mathrm{d}z $ is
finite, then the skew-symmetric family with symmetric kernel $f=g_a$
and skewing function $\Pi$ is such that Assumptions~\hyperref[as1]{(A1)}
and~\hyperref[as2plus]{(A2$^+$)}
hold, and the corresponding $2\times2$ matrix $\boldsymbol{{\Gamma
}}_{f;\varthetab_0}^{0}$ exists and is singular for any $\varthetab
_0$. If moreover $f=g_a$ also satisfies Assumption~\hyperref[as1plus]{(A1$^+$)}, then the
$3\times3$ information matrix $\boldsymbol{{\Gamma}}_{f;\varthetab
_0}$ exists, and is singular for any $\varthetab_0$. Note, however,
that the reference density for scale -- the one that, by definition,
provides the unit scale -- here is $f=g_a$.

A tale of two densities, $f$ and $g_{a_\Pi}$, is emerging, which demythifies the
seemingly singular role of the Gaussian distribution.

This treatment of the univariate case provides a good intuition for the
more complex $k$-dimensional problem where, as we shall see, the rank
of the Fisher information matrix can take any value between $k+
k(k+1)/2=k(k+3)/2$ and $2k+ k(k+1)/2=k(k+5)/2$. Since the univariate
case follows as a particular case by letting $k=1$ in the general
result of Proposition~\ref{HLmultisol} of the next section, we do not
provide a more formal statement here.
%
%s2.2 ###
\subsection{Some examples}\label{univexos}
In order to illustrate the results of the previous section, we now
apply our findings in three examples of skewing functions and determine
the exponential family with corresponding minimal sufficient statistic
and natural parameter space leading to singular Fisher information matrices.

As a first example, we propose the most usual class of skewing
functions, namely those of the form $\Pi_1(z,\delta):=\Pi(\delta z)$,
where $\Pi\dvtx\R\rightarrow[0,1]$ is a function satisfying $\Pi
(-y)+\Pi
(y)=1$ for all $y\in\R$ (hence $\Pi(0)=1/2$) and such that $\dot
{\Pi
}(0):=\mathrm{d}\Pi(y)/\mathrm{d}y\vert_{y=0}$ exists and differs from~0. Clearly, any
univariate c.d.f. could be used, in which case we retrieve the
skew-symmetric distributions of Azzalini and Capitanio \cite{AC99}, and, for
$f=\phi$ and $\Pi=\Phi$, the skew-normal distributions of Azzalini
\cite{A85}. For more examples of skewed distributions of this type, we
refer the reader to G\'{o}mez \textit{et al.} \cite{GTB07}. Straightforward
calculations show that $\psi_1(z)=\dot{\Pi}(0)z$, and hence the minimal
sufficient statistic characterizing the exponential family (\ref
{expfam}) is $\Psi_1(z)=\dot{\Pi}(0)z^2/2$. The resulting exponential
family $\mathcal{E}_{\Psi_1}$ thus is nothing but the family of
centered normal densities of the form
\[
g^{(1)}_a(z)=\exp\bigl(-a \dot{\Pi}(0)z^2/2\bigr)\bigl(2\uppi/(a\dot{\Pi}(0))\bigr)^{-1/2},
\]
with natural parameter space $\mathcal{A}_1:= \operatorname{sign}
(\dot
{\Pi}(0))\R_0^+$. Assumptions~\hyperref[as1plus]{(A1$^+$)} and \hyperref[as2plus]{(A2$^+$)} are satisfied,
hence the $3\times3$ matrix $\boldsymbol{{\Gamma}}_{f;\varthetab_0}$
exists. Thus, whenever the traditional skewing function $\Pi_1$ is
used, Gaussian kernels are the only problematic ones regarding singular
Fisher information at $\delta= 0$. This result, combined with the
popularity of $\Pi_1$ as a skewing function, explains the long-standing
belief in a particular role of the Gaussian distribution. Note that our
findings are in line with earlier ones by G\'{o}mez \textit{et al.} \cite{GTB07}, who
show that, by combining a Student kernel with $\nu$ degrees of freedom
and a skewing function of the form $\Pi_1$, Fisher information at
$\delta= 0$ is non-singular in general but becomes singular as $\nu
\rightarrow\infty$. And, more generally, our results are in total
accordance with those of Ley and Paindaveine \cite{LP10a} for the total
class of skew-symmetric distributions of this kind.

Next consider the class of skewing functions $\Pi_2(z,\delta):=\Pi
(\delta\operatorname{sign}(z)|z|^{\alpha/2}(2/\alpha)^{1/2})$
with \mbox{$\alpha>1$} and $y\mapsto\Pi(y)$ satisfying the usual conditions.
Clearly, for $\alpha=2$, $\Pi_2$ coincides with $\Pi_1$. This second
type of skewing function was used, with $\Pi=\Phi$, by Azzalini \cite{A86}
to define skew-exponential power distributions. One immediately obtains
$\psi_2(z)=\dot{\Pi}(0)\operatorname{sign}(z)|z|^{\alpha
/2}(2/\alpha
)^{1/2}$, and, consequently,
\[
\Psi_2(z)=\dot{\Pi}(0)|z|^{\alpha/2+1}(2/\alpha)^{1/2}(\alpha/2+1)^{-1}.
\]
The corresponding exponential family $\mathcal{E}_{\Psi_2}$ contains
all densities of the form
\[
g_{a}^{(2)}(z)=c\exp\bigl(-a\dot{\Pi}(0)(2/\alpha)^{1/2}(\alpha
/2+1)^{-1}|z|^{\alpha/2+1}\bigr),
\]
where $c$ is a normalization constant and $a$ again ranges over either
the positive or the negative real half line, depending on the sign of
$\dot{\Pi}(0)$. One easily can check that the complete Fisher
information matrix is well-defined in this case. DiCiccio and Monti
\cite{DM04} prove that, for $\alpha\neq2$, skew-exponential power
distributions do not suffer from singular Fisher information matrices
in the vicinity of symmetry. Our findings do not only confirm that
result, but also provide some further insight into the reasons for that
absence of singularity. Actually, the exponent of $|z|$ in $g^{(2)}_a$
has to be $\alpha/2+1$, while the symmetric kernels in skew-exponential
power distributions as defined in Azzalini \cite{A86} are of the form
$c\exp
(-|z|^\alpha/\alpha)$. Thus, while skew-normal distributions involve a
symmetric kernel and a~skewing function which are in a problematic
relationship, this is avoided with the class of skew-exponential power
distributions.

%
%f1 ###
\begin{figure}[b]

\includegraphics{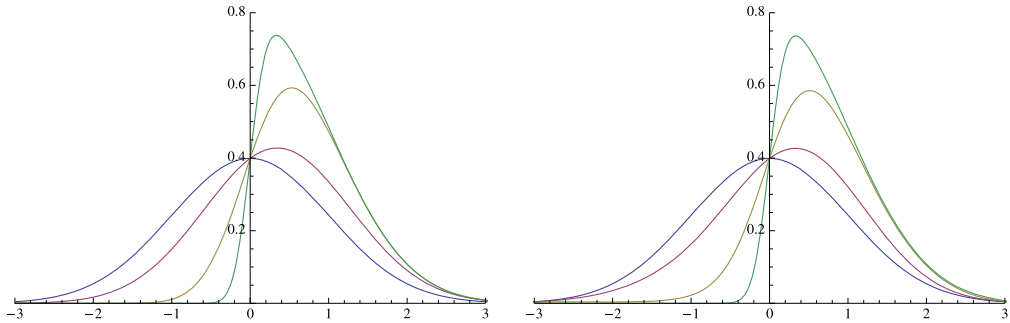}

\caption{Plots of the original Azzalini \protect\cite{A85}
skew-normal density $2\phi(x)\Phi(\delta x)$ (left) and the $\Pi
_3$-based version
$2\phi(x)\Phi(\delta\sin(x))$ (right), for $\delta=$0 (darker), 0.5,
2,
and 6 (lighter).}\label{univex}
\end{figure}

As a final example, consider skewing functions of the form $\Pi
_3(z,\delta):=\Pi(\delta\sin(z))$, with $\Pi$ belonging to the same
class of functions as in the two preceding examples. It is easy to
check that $\Pi_3$ then actually is a skewing function satisfying
Assumption~\hyperref[as2plus]{(A2$^+$)}. Direct manipulations yield $\psi_3(z)= \dot{\Pi
}(0)\sin(z)$ and $\Psi_3(z)=-\dot{\Pi}(0)\cos(z)$. The natural
parameter space $\mathcal{A}_3$ of the exponential family $\mathcal
{E}_{\Psi_3}$ corresponding to the minimal sufficient statistic $\Psi
_3$ is empty. In other words, no symmetric kernel $f$ yields a reduced
Fisher information matrix when the skewing function $\Pi_3$ is adopted.
Figure~\ref{univex} shows some of the skewed densities obtained by
combining $\Pi_3$ (for $\Pi=\Phi$) with a standard normal kernel.
Comparison with the original skew-normal distributions of Azzalini
\cite{A85} indicates that the new family, which is immune from degenerate
Fisher information problems, is nevertheless extremely close to
Azzalini's classical one.

%
%%\includegraphics[width=12.6cm]{picture1.eps}
%(yellow), and 6 (green). }%
%
%s3 ###
\section{The multivariate setup}\vspace*{3pt}
%
%s3.1 ###
\subsection{A further tale \ldots}\label{multitale}
Before starting our investigation of the multivariate case, let us
introduce some further notations required when passing from dimension 1
to $k>1$. For any given $k\times k$ matrix~$\Mb$, we denote by $\vecb
(\Mb)$ the $k^2$-vector obtained by stacking the columns of~$\Mb$ on
top of each other, and by $\vechb(\Mb)$ the $k(k+1)/2$-subvector of
$\vecb(\Mb)$ for which only upper diagonal entries in $\Mb$ are considered.
We write $\mathbf{P}_k$ for the $k(k+1)/2\times k^2$ matrix such that
$\mathbf{P}_k^\prime(\vechb\Mb)=\vecb(\Mb)$ for any symmetric
$\Mb$
and $\mathbf{I}_k$ for the $k\times k$ identity matrix.

The general multivariate skew-symmetric densities (generalizing (\ref
{skew1})) we are considering are of the form (\ref{HLSS}), with
$\varthetab$, $f$ and $\Pi$ satisfying the general conditions (a)--(c).
The symmetric kernel $f$ moreover is supposed to have identity scatter
matrix $\mathbf{I}_k$, which provides the required identification
constraint for $\Sigmab$.

As in the univariate setup, we need to impose some mild regularity
assumptions on $f$ and $\Pi$.
\begin{asum*}[(B1)]\label{asb1} The mapping ${\mathbf{z}}\mapsto
f({\mathbf{z}})$ is
differentiable, with gradient $\dot f$ such that, letting ${\boldsymbol
\varphi
}_f:=-\dot{f}/f$,
the $k\times k$ information matrix for location $\Sigmab^{-1/2}
{\boldsymbol
{\mathcal I}_f}\Sigmab^{-1/2}$, with
\[
{\boldsymbol{\mathcal I}}_f := \int_{\R^k} {\boldsymbol\varphi
}_f({\mathbf{z}}) {\boldsymbol\varphi
}_f\pr({\mathbf{z}}) f({\mathbf{z}})
\,\mathrm{d}{\mathbf{z}},
\]
is finite and invertible.
\end{asum*}
\begin{asum*}[(B1$^+$)]\label{asb1plud} Same as \hyperref[asb1]{(\textup{B}1)}, but the $k(k+1)/2\times
k(k+1)/2$ information matrix for scatter (actually, for ${\boldsymbol
\Sigma
}^{1/2}$, or more precisely, for $\operatorname{vech}(\Sigmab^{1/2} )$,
as ${\boldsymbol\Sigma}^{1/2}$ is symmetric) $\mathbf{P}_k(\Sigmab
^{-1/2}\otimes\mathbf{I}_k){{\boldsymbol{\mathcal J}}_f}(\Sigmab
^{-1/2}\otimes\mathbf{I}_k)\mathbf{P}_k^\prime$, with
\[
{{\boldsymbol{\mathcal J}}_f} := \int_{\R^k}\vecb\bigl({\mathbf
{z}}{\boldsymbol\varphi}_f^\prime
({\mathbf{z}})-\mathbf{I}_k\bigr)\bigl(\vecb\bigl({\mathbf
{z}}{\boldsymbol\varphi}_f^\prime({\mathbf{z}})-\mathbf
{I}_k\bigr)\bigr)^\prime f({\mathbf{z}}) \,\mathrm{d}{\mathbf{z}},
\]
moreover is finite and invertible.
\end{asum*}
\begin{asum*}[(B2)]\label{asb2} \emph{(i)} The mapping ${\mathbf{z}}\mapsto
\Pi({\mathbf{z}},\deltab)$
is differentiable, and has gradient $\zerob$ at $\deltab=\zerob$. \emph{(ii)}
The mapping $\deltab\mapsto\Pi({\mathbf{z}},\deltab)$
is differentiable at
$\deltab= \zerob$ for all ${\mathbf{z}}\in\R^k$, with
gradient (at $\deltab=
\zerob$) $\operatorname{grad}_\deltab\Pi({\mathbf
{z}},\deltab)|_{\deltab
=\zerob}=:\psib({\mathbf{z}})$ such that $ \psib$
admits a primitive $\Psi$,
that is, a real-valued function \mbox{${\mathbf{z}}\mapsto\Psi
({\mathbf{z}})$} such that
grad$_{{\mathbf{z}}}\Psi({\mathbf
{z}})={\boldsymbol\psi} ({\mathbf{z}})$.
\end{asum*}
\begin{asum*}[(B2$^+$)]\label{asb2plus} Same as \hyperref[asb2]{(\textup{B}2)}, but the $k\times k$ matrix
\[
\int_{\R^k}\psib({\mathbf{z}})\psib\pr(
{\mathbf{z}}) f({\mathbf{z}}) \,\mathrm{d}{\mathbf{z}}
\]
moreover is finite and invertible.
\end{asum*}

These assumptions admit the same interpretation as in the univariate
case, and basically ensure the existence of a finite Fisher information matrix.
The standardization issue also calls for the same comments as in
Section~\ref{univtale}. The interpretation of the scatter matrix~$\Sigmab$ is related to the choice of a standardization constraint on
$f$. If we impose that~$\Zb$ with p.d.f. $f$ has unit covariance matrix,
then $\Sigmab= \int_{\R^k}({\mathbf{x}}-\mub
)({\mathbf{x}}-\mub)\pr f_{\mub,\Sigmab
,\zerob}^\Pi({\mathbf{x}})\,\mathrm{d}{\mathbf{x}}$. However,
concepts of scatter that make sense irrespective of the underlying
density also can be used in this multivariate setup, such as the
celebrated \textit{Tyler matrix} ${\mathbf V}_{\mathrm{Tyler}}$ (Tyler~\cite{Ty87}),
defined as the unique symmetric positive definite matrix ${\mathbf V}$ with
$\operatorname{tr} {{\mathbf V}}=k$ satisfying
\[
\mathrm{E}\bigl[(\Xb-\mub)(\Xb-\mub)\pr/\bigl((\Xb-\mub)\pr{{\mathbf
V}}^{-1}(\Xb-\mub
)\bigr)\bigr]=k^{-1}{{\mathbf V}}.
\]
Note however that the Tyler matrix ${\mathbf V}_{\mathrm{Tyler}}$ in
fact is
a \textit{shape matrix}, not a scatter matrix: the corresponding
scatter is ${\boldsymbol\Sigma} = \sigma{\mathbf V}_{\mathrm
{Tyler}}$, with $\sigma
= k^{-1}\operatorname{tr}({\boldsymbol\Sigma})$. As in the
univariate case,
the scatter $\boldsymbol\Sigma$, for the kernel $f$, safely and
without any
loss of generality, can be fixed to identity for identification
purposes, implying that, for $f$, $\sigma$ takes value 1, while
${\mathbf
V}_{\mathrm{Tyler}}$ is an identity matrix. As in the univariate case,
this choice has no impact on the final results.

Here also, we could relax classical differentiability conditions by
considering weaker differentiability and generalized Fisher information
concepts, at the expense, however, of non-negligible technical complications.

Under Assumptions~\hyperref[asb1]{(B1)} and \hyperref[asb2]{(B2)}, the score vector ${\boldsymbol\ell
}_{f;\varthetab}$, at $\varthetab_0 := (\mub\pr, \vechb({\Sigmab}
^{1/2})\pr, \zerob\pr)^\prime$, takes the form
\begin{eqnarray*}
{\boldsymbol\ell}_{f;\varthetab_0} ({\mathbf{x}})
&:=&
\operatorname{grad}_\varthetab\log f_{\varthetab}^\Pi(
{\mathbf{x}}) \vert
_{\varthetab_0}=:
\pmatrix{
{\boldsymbol\ell}^{1\prime}_{f;\varthetab_0}({\mathbf{x}})&
{\boldsymbol\ell}^{2\prime}_{f;\varthetab_0}({\mathbf{x}})&
{\boldsymbol\ell}^{3\prime}_{f;\varthetab_0}({\mathbf{x}})
}
\pr\\
%&:=& (\begin{array}{c}
%({\operatorname{grad}}_\mub\log f_{\mub,\Sigmab,\deltab}^\Pi(\xb))|_{
%({\operatorname{grad}}_{\vechb\Sigmab^{1/2}}\log f_{\mub,\Sigmab,
%({\operatorname{grad}}_\deltab\log f_{\mub,\Sigmab,\deltab}^\Pi(
% )\vspace{1mm}\\
&=& \pmatrix{
\Sigmab^{-1/2}\boldsymbol{\varphi}_f\bigl(\Sigmab^{-1/2}(
{\mathbf{x}}-\mub)\bigr)\cr
\mathbf{P}_k(\Sigmab^{-1/2}\otimes\mathbf{I}_k)\vecb\bigl(\Sigmab
^{-1/2}({\mathbf{x}}-\mub)\boldsymbol{\varphi}\pr
_f\bigl(\Sigmab^{-1/2}({\mathbf{x}}-\mub)\bigr)-\mathbf
{I}_k\bigr)\cr
2{\boldsymbol\psi}\bigl(\Sigmab^{-1/2}({\mathbf{x}}-\mub)\bigr)
},
\end{eqnarray*}
where $\otimes$ stands for the standard Kronecker product. Note that,
for $k=1$, this score vector coincides with the one we obtained in Section~\ref
{univtale}. Under Assumptions~\hyperref[asb1plud]{(B1$^+$)} and~\hyperref[asb2plus]{(B2$^+$)}, the corresponding
Fisher information matrix
\[
\boldsymbol{{\Gamma}}_{f; \varthetab_0}:= \vert{{\Sigmab} }\vert
^{-1/2} \int_{\R^k}{\boldsymbol\ell}_{f;\varthetab_0}(
{\mathbf{x}}){\boldsymbol\ell}^{\prime
}_{f;\varthetab_0}({\mathbf{x}}) f\bigl({\Sigmab}
^{-1/2}({\mathbf{x}}- \mub) \bigr) \,\mathrm{d}{\mathbf{x}}
\]
exists and is finite, and naturally partitions into
\[
\boldsymbol{{\Gamma}}_{f; \varthetab_0}=
\pmatrix{
\boldsymbol{{\Gamma}}_{f; \varthetab_0}^{11}&\mathbf{0}&\boldsymbol
{{\Gamma}}_{f; \varthetab_0}^{13}\cr
\mathbf{0}&\boldsymbol{{\Gamma}}_{f; \varthetab_0}^{22}&\mathbf
{0}\cr
\boldsymbol{{\Gamma}}_{f; \varthetab_0}^{13\prime}&\mathbf
{0}&\boldsymbol{{\Gamma}}_{f; \varthetab_0}^{33}
},
\]
with
\begin{eqnarray*}
\boldsymbol{{\Gamma}}_{f; \varthetab_0}^{11}&=&\Sigmab^{-1/2}
{\boldsymbol
{\mathcal I}_f}\Sigmab^{-1/2},\qquad
\boldsymbol{{\Gamma}}_{f; \varthetab_0}^{22}=\mathbf
{P}_k(\Sigmab
^{-1/2}\otimes\mathbf{I}_k){{\boldsymbol{\mathcal J}}_f}(\Sigmab
^{-1/2}\otimes\mathbf{I}_k)\mathbf{P}_k^\prime, \\
\boldsymbol{{\Gamma}}_{f; \varthetab_0}^{33}&=&4\int_{\R
^k}\psib({\mathbf{z}}
)\psib\pr({\mathbf{z}}) f({\mathbf{z}})
\,\mathrm{d}{\mathbf{z}}  \quad \mbox{and} \quad  \boldsymbol{{\Gamma
}}_{f; \varthetab_0}^{13}=2\Sigmab^{-1/2} \int_{\R^k}\boldsymbol
{\varphi}_f({\mathbf{z}}
)\boldsymbol{\psi}\pr({\mathbf{z}})f(
{\mathbf{z}})\,\mathrm{d}{\mathbf{z}}.
\end{eqnarray*}

As in the univariate case, the blocks of zeroes in $\boldsymbol
{{\Gamma
}}_{f; \varthetab_0}$ readily follow from symmetry arguments and,
without loss of generality, we can focus our attention on the submatrix
\[
\boldsymbol{{\Gamma}}_{f; \varthetab_0}^{0}:= \pmatrix{
\boldsymbol{{\Gamma}}_{f; \varthetab_0}^{11}&\boldsymbol{{\Gamma}}_{f;
\varthetab_0}^{13}\cr
\boldsymbol{{\Gamma}}_{f; \varthetab_0}^{13\prime}&\boldsymbol
{{\Gamma
}}_{f; \varthetab_0}^{33}
},
\]
the existence of which only requires Assumptions~\hyperref[asb1]{(B1)} and \hyperref[asb2plus]{(B2$^+$)}. In
the univariate case, the $2\times2$ matrix $\boldsymbol{{\Gamma}}_{f;
\varthetab_0}^{0}$ was either full-rank or singular with rank 1;
here,\vspace*{2pt}
the $2k\times2k$ matrix $\boldsymbol{{\Gamma}}_{f; \varthetab_0}^{0}$
can be singular with any rank ranging from $k$ to $2k-1$ (note that the
lower bound $k$ is a direct consequence of either Assumption~\hyperref[asb1]{(B1)} or \hyperref[asb2plus]{(B2$^+$)}).

The following proposition fully characterizes, for each possible rank
$2k-m$, $m\in\{1,\ldots,k\}$, the relation between the kernel $f$ and
the skewing function $\Pi$ causing such degeneracy (for simplicity, we
restrict to a characterization of the singularity of $\boldsymbol
{{\Gamma}}_{f; \varthetab_0}^{0}$).
\begin{Prop}\label{HLmultisol}
Let the symmetric kernel $f$ and the skewing function $\Pi$ satisfy
Assumptions~\hyperref[asb1]{(\textup{B}1)} and \hyperref[asb2plus]{(\textup{B}2$^+$)}. The following statements are equivalent:
\begin{enumerate}[(ii)]
\item[(i)] the $2k\times2k$ matrix $\boldsymbol{{\Gamma}}_{f;
\varthetab_0}^{0}$ is singular with rank $2k-m$, $1\leq m\leq k$, for
any $\varthetab_0$;
\item[(ii)] denoting by $\Zb$ a random $k$-vector with p.d.f. $f$, there
exists a $k\times k$ orthogonal matrix ${\mathbf O}\pr=({\mathbf
O}_1\pr, {\mathbf
O}_2\pr)$, where ${\mathbf O}_1\pr$ and ${\mathbf O}_2\pr$ are
$k\times m$- and
$k\times(k-m)$-dimensional, respectively, such that, letting $\Yb
:={\mathbf O}\Zb$ and ${\mathbf{y}}:={\mathbf
O}{\mathbf{z}}$, for Lebesgue-almost all ${\mathbf O}
_2{\mathbf{z}}= (y _{m+1},\ldots, y _k)\pr\in\R
^{k-m}$, the density of ${\mathbf
O}_1\Zb= (Y _1,\ldots, Y _m)\pr$ conditional on ${\mathbf O}_2\Zb= (Y
_{m+1},\ldots, Y _k)\pr= (y _{m+1},\ldots,y _k)\pr$ belongs to the
exponential family
%
%e3.1 ###
\begin{eqnarray}\label{expmulti}
&&\biggl\{
(y_1,\ldots, y_m)\mapsto g_a(y_1,\ldots, y_m):=C^{-1}
%(y _{m+1},\ldots, y _k)
\exp(-a\Psi({\mathbf O}\pr{\mathbf{y}})) \Bigm|\nonumber
\\[-12pt]
\\[-6pt]
&&\quad a
\mathrm{\ such\ that\ }%=
%a(y_{m+1},\ldots, y_k)
C=C(y _{m+1},\ldots, y _k) := \int_{\R^m} \exp(
-a \Psi({\mathbf O}\pr{\mathbf{y}}))\,\mathrm{d}y_1\cdots
\,\mathrm{d}y_m<\infty
\biggr\}\qquad\nonumber
\end{eqnarray}
with parameter $a$ and minimal sufficient statistic $\Psi({\mathbf
O}\pr(Y
_1,\ldots, Y_m, y_{m+1},\ldots, y_k)\pr)$.
\end{enumerate}
\end{Prop}

Note that the natural parameter space
\[
\mathcal{A}= \mathcal{A}(y_{m+1},\ldots, y_k):=\biggl\{ a\in\R
\mbox{ such that }
\int_{\R^m} \exp(
-a \Psi({\mathbf O}\pr{\mathbf{y}}))\,\mathrm{d}y_1\cdots
\,\mathrm{d}y_m<\infty
\biggr\}
\]
of the exponential family (\ref{expmulti}) in principle also depends on
$(y_{m+1},\ldots, y_k)$. Natural parameters in exponential families
being well identified, the values $a_\Pi(y_{m+1},\ldots, y_k)$ of the
natural parameter $a$ achieving, whenever condition (ii) of
Proposition~3.1 holds, the matchings $f=g_a$, are uniquely defined for
Lebesgue-almost all $(k-m)$-tuple $(y_{m+1},\ldots, y_k)$, yielding
exponential densities $g_{a_\Pi(y_{m+1},\ldots, y_k)}$.\vspace*{\baselineskip}

Proposition~3.1 has the following straightforward corollary.
\begin{Corol}\label{HLmulticol}
\emph{(i)} Let $f$ be a symmetric kernel satisfying Assumption~\hyperref[asb1]{(\textup{B}1)}: there
exists a skewing function $\Pi_f$ such that the rank of
% the Fisher information matrix for location and asymmetry
$\boldsymbol{{\Gamma}}_{f; \varthetab_0}^{0}$ reaches its minimal
value $ k$ for any~$\varthetab_0$.\vadjust{\goodbreak}

\emph{(ii)} Let $\Pi$ be a skewing function satisfying Assumption~\hyperref[asb2]{(B2)} with
$\Psi$ such that, for some~$a_\Pi$,
\begin{enumerate}[(iib)]
\item[(iia)] ${\mathbf{z}}\mapsto g_{a_\Pi}(
{\mathbf{z}}):= C^{-1}
%(y _{m+1},\ldots, y _k)
\exp(-a_\Pi\Psi({\mathbf{z}}))$ is a p.d.f. with identity
scatter matrix, and
\item[(iib)] $\int_{\R^k}\psib({\mathbf{z}})\psib\pr
({\mathbf{z}}) f({\mathbf{z}}) \,\mathrm{d}
{\mathbf{z}}$ is
finite and invertible (meaning that \hyperref[asb2plus]{(B2$^+$)} is satisfied).
\end{enumerate}
Then, there exists a symmetric kernel $f_\Pi$ such that the rank of
% the Fisher information matrix for location and asymmetry
$\boldsymbol{{\Gamma}}_{f_\Pi; \varthetab_0}^{0}$ reaches its minimal
value $k$ for any $\varthetab_0$.
\end{Corol}
\begin{pf*}{Proof of Proposition 3.1} Clearly, $\boldsymbol{{\Gamma
}}^0_{f;\varthetab_0}$ has rank $2k-m$, $1\leq m\leq k$, if and only
if $m$ is the largest integer such that there exist $(k\times m)$
matrices ${\mathbf V}$ and ${\mathbf W}$ with $({\mathbf V}\pr,
-{\mathbf W}\pr)$ of
rank $m$ such that
%
%e3.2 ###
\begin{equation}\label{*1}
{\mathbf V}\pr{\boldsymbol\varphi}_f = {\mathbf W}\pr{\psib} \qquad\mbox
{Lebesgue-a.e.}
\end{equation}
(note that the matrix $\Sigmab^{-1/2}$ is incorporated in ${\mathbf V}$,
and hence plays no role in the characterization (\ref{*1})). Both
${\mathbf
V}$ and ${\mathbf W}$ are of maximal rank $m$. Suppose indeed that
${\mathbf
V}$ is not: then, there exists ${\boldsymbol0}\neq{\boldsymbol
\lambda}\in\R^m$ such
that ${\mathbf V}{\boldsymbol\lambda} = {\boldsymbol0}$, so that
${\boldsymbol\lambda}\pr{\mathbf
W}\pr\psib= {\boldsymbol\lambda}\pr{\mathbf V}\pr{\boldsymbol
\varphi}_f =0$
(Lebesgue-a.e.). Then, in view of Assumption~\hyperref[asb2]{(B2)}, ${\mathbf
W}{\boldsymbol\lambda
}={\boldsymbol0}$ as well, hence ${\boldsymbol\lambda}\pr({\mathbf
V}\pr,-{\mathbf W}\pr
)={\boldsymbol0}$, which contradicts the assumption that $({\mathbf
V}\pr,-{\mathbf
W}\pr)$ has rank $m$. The same reasoning holds for ${\mathbf W}$. It
follows that ${\mathbf V}$, without loss of generality, can be assumed to
be orthonormal, and therefore can be extended into an orthogonal matrix
${\mathbf O}\pr:=({\mathbf V}, {\mathbf v})$, ${\mathbf v}$ being the
$k\times(k-m)$
orthogonal complement to ${\mathbf V}$. The necessary and sufficient
condition (\ref{*1}) then takes the form
%
%e3.3 ###
\begin{equation}\label{*2}
[{\mathbf O} {\boldsymbol\varphi}_f
]_{1\ldots m} = {\mathbf W}\pr{\psib} \qquad\mbox{Lebesgue-a.e.}
\end{equation}
where $ [{\mathbf O} {\boldsymbol\varphi}_f
]_{1\ldots m} $ stands for $ {\mathbf O} {\boldsymbol\varphi}_f
$'s $m$ first rows.

Define $\Yb:={\mathbf O}\Zb$. Since $\Zb$ has density $f$, $\Yb$ has
density ${\mathbf{y}}\mapsto f^{\Yb}({\mathbf
{y}}) = f( {\mathbf O}\pr{\mathbf{y}})$. This density
$f^{\Yb}$ has gradient $\dot{f}^\Yb$ and score $\varphib_{f^{\Yb
}}$, with
\[
\varphib_{f^{\Yb}}({\mathbf{y}}):=-\dot{f}^\Yb
({\mathbf{y}}) / f^\Yb({\mathbf{y}}) =
-{\mathbf
O} \dot{f}({\mathbf O}\pr{\mathbf{y}})/{f}({\mathbf
O}\pr{\mathbf{y}}) = {\mathbf O} \varphib_f
({\mathbf O}\pr{\mathbf{y}}).
\]
This, combined with (\ref{*2}), yields
\[
[
\varphib_{f^{\Yb}}({\mathbf{y}})
]_{1\ldots m} = {\mathbf W}\pr{\psib}({\mathbf O}\pr
{\mathbf{y}}) \qquad\mbox{Lebesgue-a.e.}
\]
or, more explicitly,
%
%e3.4 ###
\begin{equation}\label{*3}
\pmatrix{
\partial_{y_1}\log f^\Yb({\mathbf{y}})\cr
\vdots\cr
\partial_{y_m}\log f^\Yb({\mathbf{y}})
}
=
-{\mathbf W}\pr{\psib}({\mathbf O}\pr{\mathbf{y}})
\qquad\mbox{Lebesgue-a.e.}
\end{equation}
As a function of $(y_1,\ldots, y_m)$, the left-hand side in (\ref{*3})
has primitive
\[
\log f^\Yb(y_1,\ldots,y_m, y_{m+1},\allowbreak\ldots, y_k) +
c(y_{m+1},\ldots, y_k),
\]
where the ``integration constant'' $c$ is an
arbitrary function of $(y_{m+1},\ldots, y_k)$. The right-hand side
therefore has the same primitive,\vadjust{\goodbreak} still up to an additive
$c(y_{m+1},\ldots, y_k)$. Now, partitioning ${\mathbf O}\pr$ into
$({\mathbf
O}\pr_1 , {\mathbf O}\pr_2)$ where ${\mathbf O}\pr_1$ and ${\mathbf
O}\pr_2$ are
$k\times m$ and $k\times(k-m)$, respectively, a necessary condition for
\[
(y_{1},\ldots, y_m)\mapsto{\mathbf W}\pr\psib\bigl({\mathbf O}\pr
_1(y_{1},\ldots
, y_m)\pr+ {\mathbf O}\pr_2(y_{m+1},\ldots, y_k)\pr\bigr)
\]
to be the gradient of a scalar function is ${\mathbf W}\pr= a{\mathbf O}_1$
for some $a =a(y_{m+1}, \ldots, y_k) \in\R$: in view of Assumption~\hyperref[asb2]{(B2)},
a primitive of
\[
(y_{1},\ldots, y_m)\mapsto a{\mathbf O}_1\psib\bigl({\mathbf O}\pr
_1(y_{1},\ldots
, y_m)\pr+ {\mathbf O}\pr_2(y_{m+1},\ldots, y_k)\pr\bigr)
\]
is then $a\Psi({\mathbf O}\pr_1(y_{1},\ldots, y_m)\pr+ {\mathbf
O}\pr
_2(y_{m+1},\ldots, y_k)\pr)$, up to the usual additive
constant -- here, an arbitrary function of $(y_{m+1}, \ldots, y_k)$.
The necessary and sufficient condition (\ref{*3}) thus takes the
further form
\[
f^\Yb({\mathbf{y}})
=
\exp(
-c(y_{m+1},\ldots, y_k)
)
\exp\bigl(
-a \Psi\bigl({\mathbf O}\pr_1(y_{1},\ldots, y_m)\pr+ {\mathbf O}\pr
_2(y_{m+1},\ldots, y_k)\pr\bigr)
\bigr)
\]
for some $a = a(y_{m+1},\ldots, y_k)\in\R$; in other words, the
conditional density of $(Y_1,\ldots,Y_m)\pr$ given $(Y_{m+1},\ldots
,Y_k)\pr=(y_{m+1},\ldots,y_k)\pr$ is
%e3.5 ###
\begin{eqnarray}\label{match}
&&f^{(Y_1,\ldots,Y_m)\pr\mid (Y_{m+1},\ldots,Y_k)\pr=(y_{m+1},\ldots
,y_k)\pr}(y_1,\ldots,y_m)\nonumber\\
&&\quad=f^\Yb(y_1,\ldots,y_m,y_{m+1},\ldots, y_k ) \Big/\int_{\R^m} f^\Yb
(y_1,\ldots,y_m,y_{m+1},\ldots, y_k) \,\mathrm{d}y_1\cdots \,\mathrm{d}y_m\hspace*{12pt}\\
&&\quad=
C(y_{m+1},\ldots, y_k)
\exp\bigl(
-a
%(y_{m+1},\ldots, y_k)
\Psi\bigl({\mathbf O}\pr_1(y_{1},\ldots, y_m)\pr+ {\mathbf O}\pr
_2(y_{m+1},\ldots, y_k)\pr\bigr)
\bigr),\nonumber
\end{eqnarray}
where $C^{-1}(y_{m+1},\ldots, y_k):=\hspace*{-1pt}
\int_{\R^m}\exp(
-a \Psi({\mathbf O}\pr_1(y_{1},\ldots, y_m)\pr+ {\mathbf O}\pr
_2(y_{m+1},\ldots,  y_k)\pr)
)\,\mathrm{d}y_1\cdots \,\mathrm{d}y_m$, for some $a=a(y_{m+1},\ldots, y_k)\in\R$.

Summing up, there exists an orthogonal matrix ${\mathbf O}\pr=
({\mathbf O}\pr
_1 , {\mathbf O}\pr_2)$ such that, for any $(y_{m+1}, \ldots,\allowbreak y_k)\pr
\in\R
^{k-m}$, the density of ${\mathbf O}_1\Zb=:(Y_1,\ldots, Y_m)\pr$
conditional on ${\mathbf O}_2\Zb= (y_{m+1}, \ldots, y_k)\pr$ belongs to
the exponential family with minimal sufficient statistic
\[
\Psi\bigl({\mathbf O}\pr_1(Y_{1},\ldots, Y_m)\pr+ {\mathbf O}\pr
_2(y_{m+1},\ldots
, y_k)\pr\bigr)
,
\]
as was to be proved.
\end{pf*}

So far, we have formally solved the singularity problem for the
$2k\times2k$ information matrix $\boldsymbol{{\Gamma}}_{f;
\varthetab
_0}^{0}$. As in the univariate case, the singularity problem for the
full $k(k+5)/2\times k(k+5)/2$ information matrix $\boldsymbol{{\Gamma
}}_{f; \varthetab_0}$ is slightly different. Indeed, the existence of
$\boldsymbol{{\Gamma}}_{f; \varthetab_0}$ requires the stronger
Assumption~\hyperref[asb1plud]{(B1$^+$)}, as the information for scatter, which is not
present in $\boldsymbol{{\Gamma}}_{f; \varthetab_0}^{0}$, has\vspace*{2pt} to exist
as well; this adds a further condition on the exponential family in
Proposition~3.1. Nevertheless, there is no fundamental difference
between the two setups: it only could happen that a solution to the
singularity problem of $\boldsymbol{{\Gamma}}_{f; \varthetab_0}^{0}$ is
not a solution of the larger problem because the matrix $\boldsymbol
{{\Gamma}}_{f; \varthetab_0}$ simply does not exist, hence cannot be
singular. This explains why, for the sake of simplicity, we state the
results of this section in terms of $\boldsymbol{{\Gamma}}_{f;
\varthetab_0}^{0}$. The message is clear: the tale of two densities has
turned into a more elaborate plot, starring a much larger number of actors.
%
%s3.2 ###
\subsection{Further examples}\label{multiexos}
As in the univariate case, we now analyze three concrete examples of
skewing functions in the light of the findings of the previous section,
which provides the theoretical statement in Proposition~\ref
{HLmultisol} with some further intuition.

The first example is the natural extension of the univariate skewing
function $\Pi_1$ to the multivariate context, with $\Pi
_1^{(k)}({\mathbf{z}}
,\deltab):=\Pi(\deltab\pr{\mathbf{z}})$, where $\Pi
\dvtx\R\rightarrow[0,1]$
satisfies exactly the same conditions as in Section~\ref{univexos}. The
resulting class of skewing functions $\Pi_1^{(k)}$ is the most common
one in the literature. A skewing function $\Pi=\Phi$ combined with a
multinormal kernel $f=\phi_k$ yields the class of skew-multinormal
densities of Azzalini and Dalla Valle~\cite{AD96}. When $f$ is only required
to be spherically symmetric and the skewing function~$\Pi$ is a
univariate symmetric c.d.f., we obtain the class of skew-elliptical
distributions as defined by Azzalini and Capitanio \cite{AC99}, itself a
subclass of the generalized skew-elliptical distributions of Genton and
Loperfido \cite{GL05} where $\Pi$ is left unspecified. Finally, relaxing the
assumption of spherical symmetry into the weaker assumption of central
symmetry, we retrieve the popular class of multivariate skew-symmetric
distributions analyzed in Ley and Paindaveine \cite{LP10a}.

Direct calculation yields $\psib^{(k)}_1({\mathbf
{z}})=\dot{\Pi}(0){\mathbf{z}}$, hence,
writing ${\mathbf{z}}=({\mathbf{z}}_1\pr
,{\mathbf{z}}_2\pr)\pr$ with ${\mathbf
{z}}_1\in\R^m$ and \mbox{${\mathbf{z}}_2\in\R
^{k-m}$}, $m=1,\ldots, k$, we obtain minimal sufficient statistics of
the form
\[
\Psi^{(k)}_1({\mathbf O}\pr(\Zb_1\pr, {\mathbf
{z}}_2\pr)\pr)=
\dot{\Pi}(0)(\Zb_1\pr\Zb_1/2+{\mathbf{z}}_2\pr
{\mathbf{z}}_2/2)
\]
for a $k\times k$ orthogonal matrix decomposing into ${\mathbf O}\pr
=({\mathbf
O}_1\pr,{\mathbf O}_2\pr)$. Quite nicely, the possibility of
separating the
vectors $\Zb_1$ and ${\mathbf{z}}_2$ in $\Psi
^{(k)}_1({\mathbf O}\pr(\Zb_1\pr, {\mathbf{z}}
_2\pr)\pr)$
allows us to express the corresponding exponential densities in terms
of ${\mathbf{z}}_1$ only, yielding the $m$-dimensional
Gaussian densities
%The resulting exponential densities are $m$-dimensional Gaussian, of
%the form
%
\[
{\mathbf{z}}_1\mapsto\exp\bigl(-a\dot{\Pi}(0)
{\mathbf{z}}_1\pr{\mathbf{z}}_1/2\bigr)\bigl(2\uppi/(a\dot{\Pi
}(0))\bigr)^{-m/2}.
\]
As in the univariate case, the sign of $a$ is the same as that of $\dot
{\Pi}(0)$. Degenerate information thus takes place iff, for some
adequate rotation $\Ob\Zb$ of $\Zb\sim f$, the $m$-dimensional marginal
distribution of $[\Ob\Zb]_{1\ldots m}$ is standard $m$-variate normal.
Note that this does not imply $k$-variate normal distributions.
Consider, for example, a random $k$-vector whose first~$m$ components
are i.i.d. standard Gaussian, and independent of the remaining
$k-m$ ones, themselves i.i.d. with some other standardized
univariate symmetric distribution. In such a case, the conditional
distribution of the~$m$ first components given the~\mbox{$k-m$} last ones
belongs to the exponential family of distributions just described.
Thus, contrary to the univariate setup, multinormal densities are not
the only symmetric kernels leading to singular Fisher information when
combined with the skewing functions $\Pi^{(k)}_1$. Multinormal kernels,
however, are the only ones for which Fisher information has minimal
rank (corresponding to $m=k$). All this is in total accordance with
earlier findings by Ley and Paindaveine \cite{LP10a}, who examine in detail
the singularity issues related to skew-symmetric distributions
generated via $\Pi^{(k)}_1$. We therefore refer the reader to that
reference for more details about the skewing functions $\Pi^{(k)}_1$,
especially so for the special case of skew-elliptical
distributions.\vadjust{\goodbreak}

Our second example corresponds to another classical type of skewing
functions, namely
%
%e3.6 ###
\begin{equation}\label{ex2}
\Pi^{(k)}_2({\mathbf{z}},\deltab):=\Pi\bigl(\deltab\pr
{\mathbf{z}}(\nu+k)^{1/2}({\mathbf{z}}\pr
{\mathbf{z}}+\nu
)^{-1/2} \bigr),
\end{equation}
where $\Pi$ satisfies the same properties as above, and $\nu>0$.
Clearly, as $\nu\rightarrow\infty$, the skewing functions $\Pi^{(k)}_2$
tend to skewing functions of the $\Pi^{(k)}_1$ type just considered.
When $\Pi$ in~(\ref{ex2}) corresponds to the c.d.f. $T_1( \cdot,\nu+k)$
of a Student variable with $\nu+k$ degrees of freedom, and the
symmetric kernel used is a $k$-dimensional $t$ variable with $\nu$
degrees of freedom, then we obtain the celebrated multivariate skew-$t$
distributions of Azzalini and Capitanio~\cite{AC03} -- up to some minor
details, since their non-standardized skewing functions are of the form
\[
T_1\bigl(\deltab\pr\omegab^{-1}({\mathbf{x}}-\mub)(\nu
+k)^{1/2}\bigl(({\mathbf{x}}-\mub)\pr\Sigmab
^{-1}({\mathbf{x}}-\mub)+\nu\bigr)^{-1/2};\nu+k\bigr),
\]
with $\omegab=\operatorname{diag}(\Sigmab_{11},\ldots,\Sigmab
_{kk})^{1/2}$. Elementary calculation yields
\[
\psib^{(k)}_2({\mathbf{z}})=\dot{\Pi}(0)
{\mathbf{z}}(\nu+k)^{1/2}({\mathbf{z}}\pr
{\mathbf{z}}+\nu)^{-1/2},
\]
hence minimal sufficient statistics and exponential densities of the
form
\[
\Psi^{(k)}_2=\dot{\Pi}(0)(\nu+k)^{1/2}({\mathbf{z}}\pr
{\mathbf{z}}+\nu)^{1/2}
\]
and
%
%e3.7 ###
\begin{eqnarray}\label{multit}
&&\exp\bigl(-a\dot{\Pi}(0)(\nu+k)^{1/2}({\mathbf{z}}\pr
{\mathbf{z}}+\nu)^{1/2}\bigr) \nonumber
\\[-8pt]
\\[-8pt]
&&\quad{}\Big/ \int_{\R
^m}\exp\bigl(-a\dot{\Pi}(0)(\nu+k)^{1/2}({\mathbf{z}}\pr
{\mathbf{z}}+\nu)^{1/2}\bigr)\,\mathrm{d}z_1\cdots \,\mathrm{d}z_m,\nonumber
\end{eqnarray}
respectively. Here again, the sign of $a$ is determined by the sign of
$\dot{\Pi}(0)$. Azzalini and Genton \cite{AG08} conjecture that, as long as
$\nu$ is finite, multivariate skew-$t$ distributions should be free of
singularity problems. DiCiccio and Monti \cite{DM10} prove the conjecture in
the univariate case, Ley and Paindaveine \cite{LP10b} in any dimension $k$.
Proposition~3.1 confirms those earlier results, as (\ref{multit}),
whatever the value of $a$, cannot be derived from a $k$-dimensional $t$
distribution with $\nu$ degrees of freedom. Actually, letting $\Xb
=(\Xb
_1\pr,\Xb_2\pr)\pr$ follow a $k$-variate $t$ distribution where
$\Xb_1$
and $\Xb_2$, respectively, are $m$- and $(k-m)$-dimensional random
vectors, it can be shown that the density of $\Xb_1\vert \Xb
_2={\mathbf{x}}_2$
cannot be of the form (\ref{multit}).

We conclude this section with a possible extension of the
singularity-free univariate skewing function $\Pi_3$ of Section~\ref{univexos}.
Consider $\Pi^{(k)}_3({\mathbf{z}},\deltab):=\Pi
(\deltab\pr\operatorname{Sin}({\mathbf{z}}
))$, with $\Pi$ defined as above and $\operatorname{Sin}(
{\mathbf{z}}):=(\sin
(z_1),\ldots,\sin(z_k))\pr$. Checking the validity of Assumption~\hyperref[asb2plus]{(B2$^+$)} is immediate, and one also directly obtains that $\psib
^{(k)}_3({\mathbf{z}})=\dot{\Pi}(0)\operatorname{Sin}
({\mathbf{z}})$ and $\Psi^{(k)}_3
=-\dot{\Pi}(0)(\cos(z_1)+\cdots+\cos(z_k))$. The same reasoning as for
$\Pi_3$ readily yields that the natural parameter space related to the
exponential family with minimal sufficient statistic $\Psi^{(k)}_3$ is
empty, hence skewing functions of the type $\Pi^{(k)}_3$ can be used
without worrying about possibly singular Fisher information.
%
%s4 ###
\section{Final comments}\label{fcomments}
In this paper, we fully dispel the widespread opinion that Gaussian
densities, in the context of skew-symmetric distributions, constitute
an intriguing worst-case situation, being the only ones (possibly,
after restriction to linear subspaces) leading to degenerate Fisher
information matrices in the vicinity of symmetry.
%We have completely solved in the previous sections the long-standing
%problem of a possible singularity of the Fisher information matrix, in
%the vicinity of symmetry, within the general class of multivariate
%skew-symmetric distributions.
Our main result provides a complete characterization of that
information degeneracy phenomenon, which generalizes and extends all
previous results of that type, and highlights the link between the
symmetric kernel and the skewing function causing singularity. We also
show how that link, in the univariate as well as in the multivariate
case, can be described as a mismatch between two densities, in which
the Gaussian distribution plays no particular role. By avoiding such
mismatch, one can deal with skew-symmetric distributions without
worrying about singular Fisher information and its consequences.
\section*{Acknowledgements}
Marc Hallin is also member of the Acad\' emie Royale de Belgique and
ECORE, and an extra-muros Fellow of CentER, Tilburg University.
His research is supported by the Sonderforschungsbereich ``Statistical
modeling of nonlinear dynamic processes'' (SFB 823) of the German
Research Foundation (Deutsche Forschungsgemeinschaft) and a Discovery
Grant of the Australian Research Council. The financial support and
hospitality of ORFE and the Bendheim Center at Princeton University,
where part of this work was
completed, is gratefully acknowledged. Christophe Ley thanks the Fonds
National de la Recherche Scientifique, Communaut\'{e} fran\c{c}aise de
Belgique, for support via a Mandat d'Aspirant FNRS. Both authors would
like to thank two anonymous referees for helpful comments that led to
an improvement of the paper.

% imsref loaded by audrone.aklyte, 2012-03-16 13:39:00
% imsref loaded by audrone.aklyte, 2012-03-16 14:08:10

\printhistory

\end{document}